\documentclass[12pt]{birkjour}
\usepackage{amsfonts}
\usepackage{graphicx,subfigure}
\usepackage{color}
\usepackage[dvipsnames]{xcolor}
\theoremstyle{plain}
\newtheorem{thm}{Theorem}[section]
\newtheorem{lem}[thm]{Lemma}

\newtheorem{que}[thm]{Question}

\theoremstyle{definition}

\newtheorem{remark}[thm]{Remark}
\begin{document}
\title[Simply Connected 3-Manifolds with  Ends of Specified Genus]
{Simply Connected 3-Manifolds with a Dense Set of Ends of Specified Genus}
\author{Dennis J. Garity}
\address{Mathematics Department, Oregon State University\\
Corvallis, OR 97331, U.S.A.}
\email{garity@math.oregonstate.edu}
\author{Du\v{s}an D. Repov\v{s}}
\address{Faculty of Education, Faculty of  Mathematics and Physics,
University of Ljubljana,  
Ljubljana, SI-1000, Slovenia}
\email{dusan.repovs@guest.arnes.si}
\thanks{The first author was supported in part by NSF grant DMS 0139678. The second author was supported in part by SRA grant P1-0292.
Both authors were supported in part by SRA grants J1-7025, J1-6721, and BI-US/15-16-029. We would like to thank the referee for remarks and suggestions.
}
\date{\today}
\begin{abstract}
We show that for every sequence  $(n_i)$, where  each $n_i$ is either an integer greater than 1 or is $\infty$,
there exists a simply connected open 3-manifold $M$ with a countable dense set of ends $\{e_i\}$
 so that, for every $i$,  the genus of end $e_i$ is equal to $n_i$. In addition, the genus of the ends not in the dense set is shown to be less than or equal to 2. These simply connected 3-manifolds are constructed as the complements of certain Cantor sets in $S^3$. The methods used require careful analysis of the genera of ends and new techniques for dealing with infinite genus. 
 \end{abstract}
\maketitle
\section{Introduction}
\label{introsec}
The \emph{ends} of an open 3-manifold $M$ are determined by any properly nested  sequence
$C_1\subset C_2\subset \ldots$ of compact 3-manifolds  $C_i$ whose union is $M$.
Such a sequence is called an \emph{exhaustion} of $M$.
The components of $M\setminus C_i$ form an inverse sequence with the bonding maps given by inclusions.

Each end $e$ is a point in this inverse sequence. An end $e$ is thus associated with a sequence $D_1\supset D_2\ldots$ where  $D_i$ is a component of $M\setminus C_i$. The \emph{genus} of the end $e$ associated with the sequence $(C_i)$, denoted by $g(e,(C_i)),$ is defined to be  $\sup\{g(D_i)\}$, where $g(D_i)$ is the genus of the boundary of  $D_i$. 
The \emph{genus} of the end $e$, denoted by $g(e)$, is the minimum of $g(e,(C_i))$, taken over all possible exhaustions $(C_i)$ of the manifold $M$.

The \emph{endpoint (or Freudenthal) compactification} of $M$ adds a point to $M$ for each end.  We say that  \emph{a set of ends is dense} (in the set of all ends) if it is dense in the set of all ends of the endpoint compactification.  For background on  endpoint compactifications and the theory of
ends, see \cite{Di68}, \cite{Fr42},  \cite{Gu16}, and \cite{Si65}. 
 For an
alternate proof using defining sequences of the result that every
homeomorphism of the open $3$-manifold extends to a homeomorphism
of its endpoint compactification, see \cite{GR13}.

It is easy to construct 3-manifolds with a finite number of ends of  specified genus. For example, if $H$ is a surface of genus $k$, then $H\times \mathbb{R}$ has two ends of genus $k$. It is more difficult to construct simply connected examples or examples with more complicated end structures. The examples we produce will be complements of certain Cantor sets in $S^3$. The Cantor sets will have simply connected complements and the end sets of the complements will be these Cantor sets. See Kirkor \cite{Ki58}, DeGryse and Osborne \cite{DO74},
Ancel and Starbird \cite{AS89}, and Wright \cite{Wr89} for
further discussion of wild (non-standard) Cantor sets with simply connected
complement.
For more recent results on Cantor sets with simply connected complements, see \cite{GRWZ11} and  \cite{GRW14}.

The main theorem we will prove in this paper is the following.

\begin{thm}
\label{EndTheorem}Let $\mathfrak{S}=(n_1,n_2, \ldots)$ be a sequence where each $n_i$ is either an integer greater than  $1$ or is $\infty$. Then there is a simply connected open 3-manifold $M=M_{\mathfrak{S}}$ with 
uncountably many ends so that the following holds: There is a countable dense set $D$ of ends 
$\{e_1,e_2, \ldots\}$ in the endpoint compactification of $M$ so that, for every $i$, the genus of $e_i$ is $n_i$. The genus of each end not in $D$ is less than or equal to 2. 
\end{thm}

This theorem will follow as a corollary of the following result about Cantor sets in $S^3$.

\begin{thm}
\label{CantorTheorem}Let $\mathfrak{S}=(n_1,n_2, \ldots)$ be a sequence where each $n_i$ is either an integer greater than  $1$ or is $\infty$. Then there is a Cantor set  $C=C_{\mathfrak{S}}$ 
in $S^3$ with simply connected complement  so that the following holds: There is a countable dense set $D$ of points
$\{x_1,x_2, \ldots\}$ in $C$ so that, for every $i$,  the local genus of $C$ at $x_i$ is $n_i$. The local genus of $C$ at each point  not in $D$ is less than or equal to 2. 
\end{thm}

\emph{Section \ref{localgenussec}} provides details about the local genus of points in a Cantor set in $S^3$
(introduced in \cite{Ze05}) and relates this to the genus of an end of the 3-manifold that is the complement of the Cantor set. \emph{Section \ref{replacementsec}} gives two replacement constructions we will need in the construction of our examples. \emph{Section \ref{constructionsec}} gives the main construction of our examples. \emph{Section \ref{proofsec}} shows that the examples we construct give a proof of the theorems from \emph{Section \ref{introsec}}. \emph{Section \ref{questionsec}} lists some question arising from our results.

 For the remainder of this paper, we assume all surfaces referred to are closed (i.e. compact, connected and without boundary) and that all 1-handles are orientable and unknotted.

\section{Local genus of points in a Cantor set}
\label{localgenussec}

We review the definition and some facts from
\cite{Ze05} about the local genus of a point in a Cantor set  in $S^3$. At the end of this section, we relate the local genus of a point $x$ in a Cantor set $C\subset S^3$ to the genus of an end $x$ of the complementary 3-manifold $S^3 \setminus C$. We use $\rm Int X$ and $\rm Fr X$ to denote the topological interior and boundary of a subset $X$ of a space $Y$.

Let $C$ be a Cantor set in $S^3$. A \emph{defining sequence} for $C$ is a
nested sequence $(M_i)$ of compacta  $M_i\subset C$ whose intersection is $C$ such that each $M_{i}$ consists of pairwise disjoint handlebodies and so, for every $i$,  that  $M_{i+1}\subset \mathop{\rm Int}\nolimits M_{i}$. The \emph{genus} of a handlebody H is denoted by $g({H})$.

Let ${ \mathcal{D}}(C)$ be the set of all defining sequences for  $C$.
Let $(M_i)\in{\mathcal{D}}(C)$. For any  $x\in C$ we denote
by $M_i^{\{x\}}$ the  component of $M_i$ which contains
$x$. Define
\begin{center}
$
g_x(C;(M_i)) = \sup\{g(M_i^{\{x\}});\ i\geq0\}\ \ \mbox{ and}$ \\
$g_x(C) = \inf\{ g_x(C;(M_i));\ (M_i) \in {\mathcal{D}}(C)\}
$.
\end{center}

The number $g_x(C)$ is called   \emph{the local genus of the Cantor set
$C$ at the point $x$}.

\begin{remark} 
\label{SubCantorSet}
This definition immediately implies that $g_x(D)\leq g_x(C)$ if $D$ is a Cantor set contained in $C$.
\end{remark}

Determining $g_x(C)$ using the
definition can be difficult. If  a defining
sequence for $C$ is given, one can easily determine an upper bound. The idea of
slicing discs introduced in \cite{Ba92} can be used to derive
the following addition theorem for local genus. This can then be
used for establishing the exact local genus. See \cite[Theorem
14]{Ze05} for details.

\begin{thm}
\label{Slicing} Let $X,Y\subset S^3$ be Cantor sets and $p$ a
point in  $X\cap Y$.  Suppose there exist a 3-ball $B$ and a 2-disc $D\subset B $ such that

\begin{enumerate}
\item $p\in\mathop{\rm Int}\nolimits B$, $\mathop{\rm
Fr}\nolimits D=D\cap\mathop{\rm Fr}\nolimits B$,
$D\cap (X\cup Y)=\{p\}$; and

\item $X\cap B\subset B_X\cup\{p\}$ and $Y\cap B\subset
B_Y\cup\{p\}$ where $B_X$ and $B_Y$ are the components of
$B\setminus D$.
\end{enumerate}

Then $g_p(X\cup Y)=g_p(X)+g_p(Y)$.
\end{thm}

The 2-disc $D$ in the above theorem is called a \emph{slicing disc} for
the Cantor set $X \cup Y$.

The lemma below points out the relationship between genus of an end and local genus of points in a Cantor set in the context we are discussing.

\begin{lem}
\label{equivalence}
Let $C\subset S^3$ be a Cantor set. Let $M=S^3\setminus C$. Then the endpoint compactification of $M$ is $S^3$ with the ends of $M$ corresponding to the points in $C$. The genus $g(e)$ of an end $e$ of $M$ is equal to the local genus $g_e(C)$ of $e$ in $C$.
\end{lem}

The proof of this lemma is just an exercise in the definitions. This follows since every defining sequence $(M_i)$ for $C$ corresponds to an exhaustion $C_1\subset C_2\subset \ldots$ of $M$ where
$C_i=\overline{S^3\setminus M_i}$.

\section{Replacement Constructions}
\label{replacementsec}
In this section, given a solid handlebody $H\subset S^3$, we define two possible constructions that replace $H$ by one or more handlebodies contained in  the interior of  $H$. These replacement constructions will be used in the next section in constructing a defining sequence for a certain Cantor set with the properties needed for the proof of the main theorem.

\subsection{Genus Replacement} 
\begin{figure}[ht]
\begin{center}
  \subfigure[Original Handlebody]%
    {%
    \label{Replacement1}
    \includegraphics[width=0.28\textwidth]{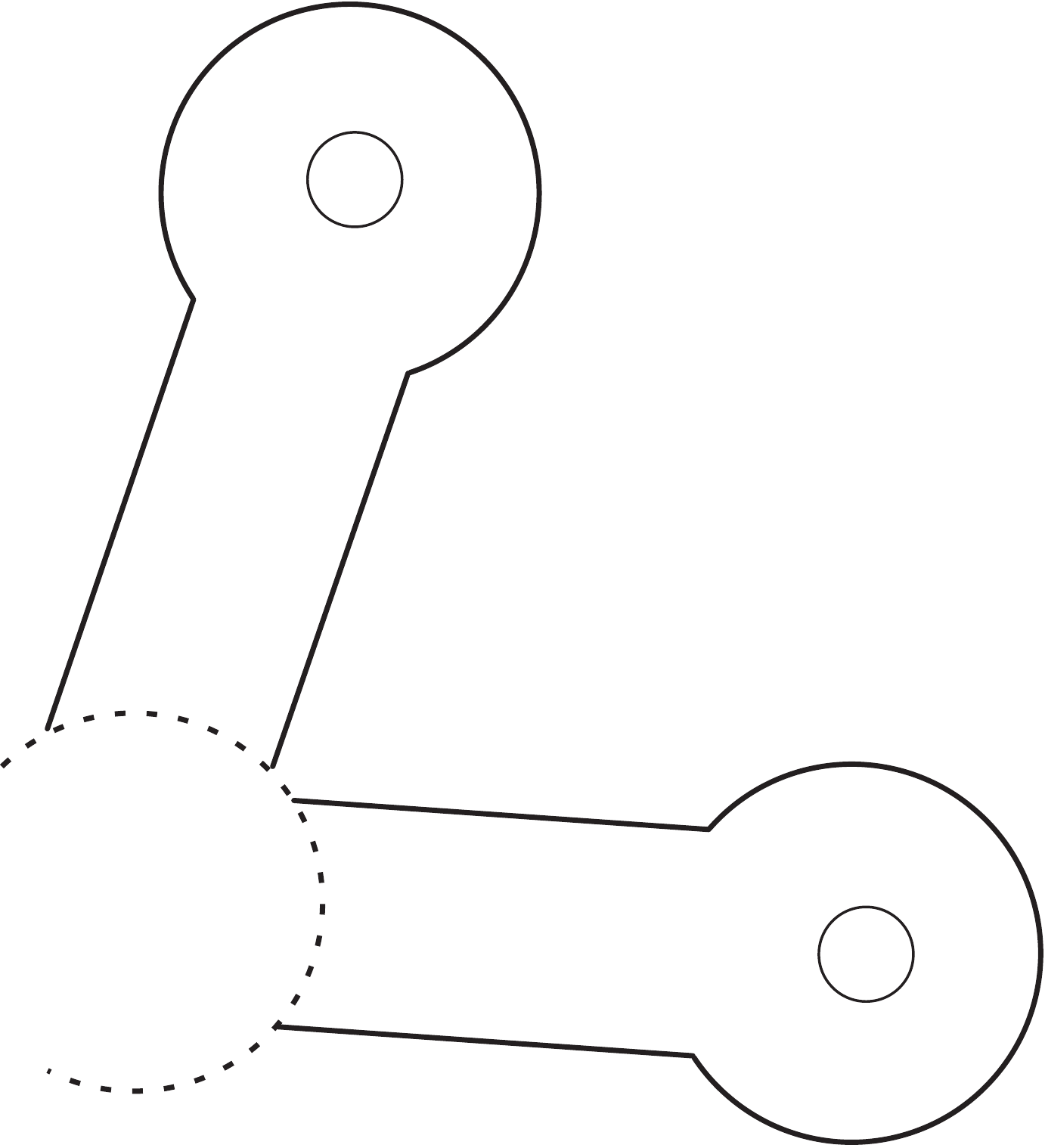}
    }%
  \quad
  \subfigure[Original and Replacement]%
    {%
    \label{Replacement2}
    \includegraphics[width=0.28\textwidth]{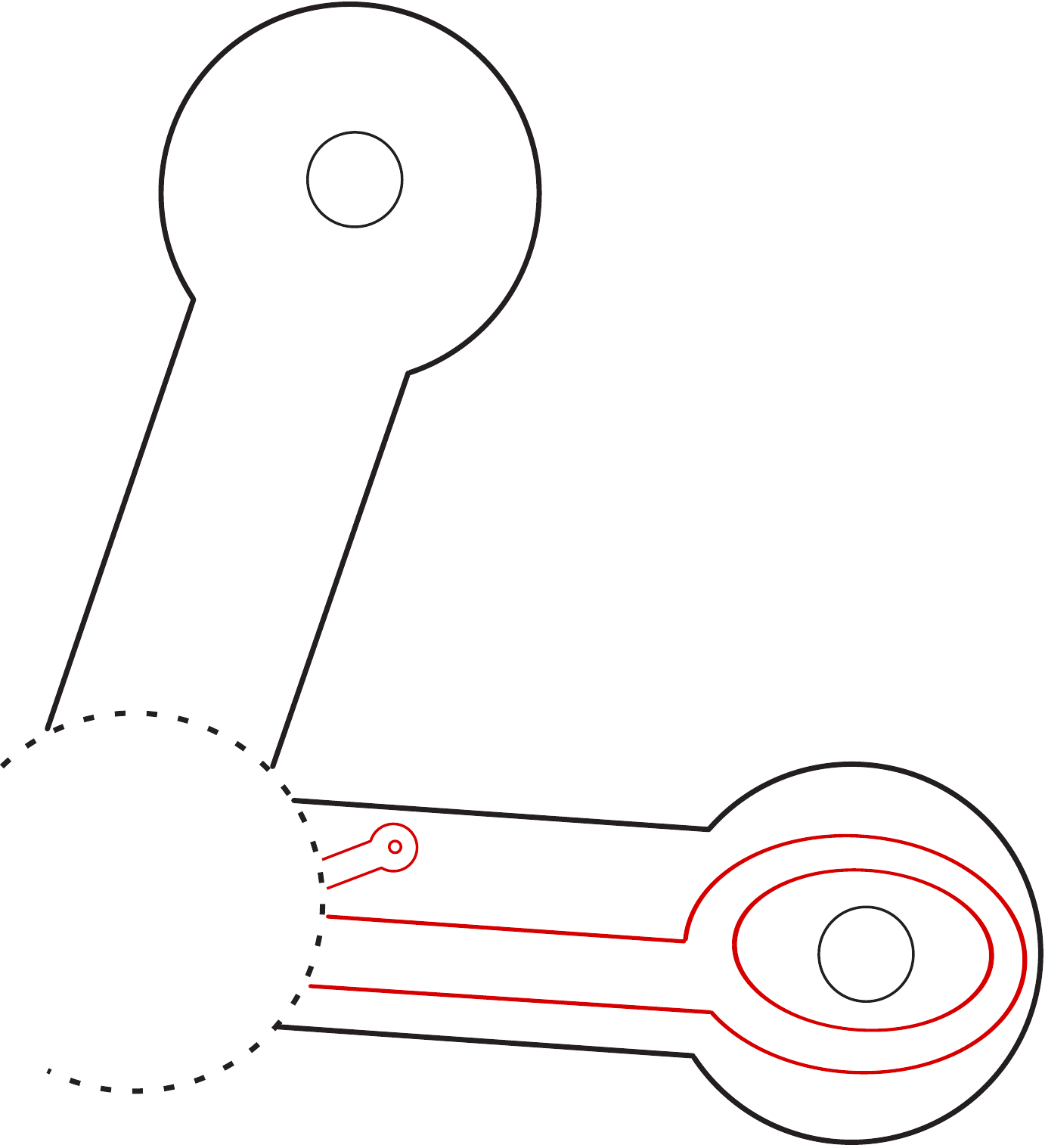}
    }
    \quad
    \subfigure[Replacement Handlebody]%
    {%
    \label{Replacement3}
    \includegraphics[width=0.28\textwidth]{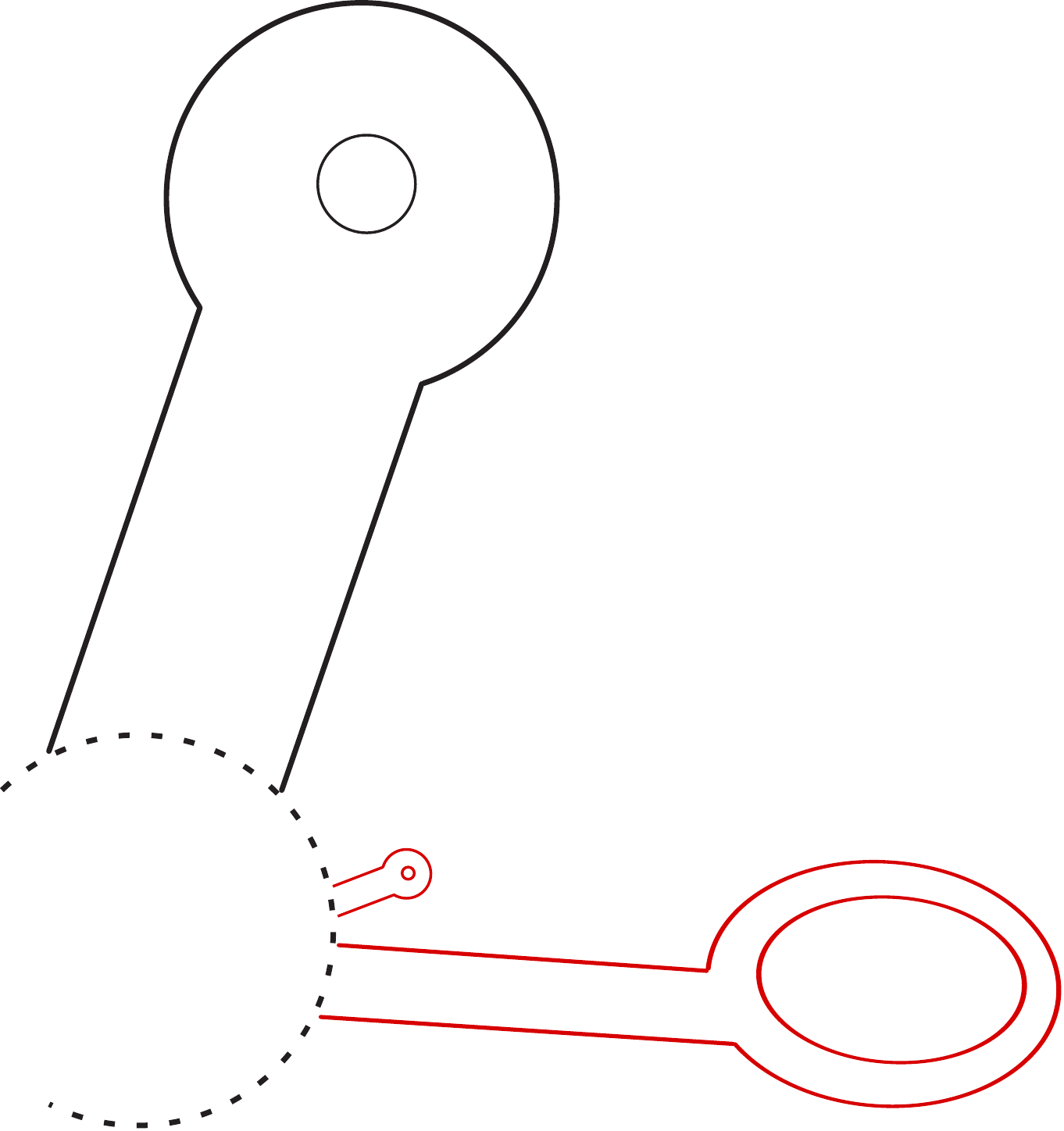}
}%
  \end{center}
  \caption{%
    Genus Replacement}%
  \label{GenusR}
\end{figure}

(See Figure \ref{GenusR}.) Let $H$ be a unknotted genus $k$ solid handlebody in $S^3$. View $H$ as a solid ball with $k$ 1-handles attached as in Figure \ref{Replacement1}. 

The genus replacement construction replaces $H$ by a handlebody $H^{\prime}\subset H$ of genus $k+1$.  One of the handles of $H$ is replaced by two smaller handles. One of the new handles is a parallel interior copy of the replaced handle. The other new handle is contractible in $H$ and lies in the interior of the replaced handle as in Figure \ref{GenusR}. As a final step, remove a small neighborhood of the boundary of the new handlebody to guarantee that $H^{\prime}$ is in the interior of $H$ as needed for a defining sequence.

\subsection{Size Replacement} 
\begin{figure}[ht]
\begin{center}
  \subfigure[General Pattern]%
    {%
    \label{SizeR1}
    \includegraphics[width=0.52\textwidth]{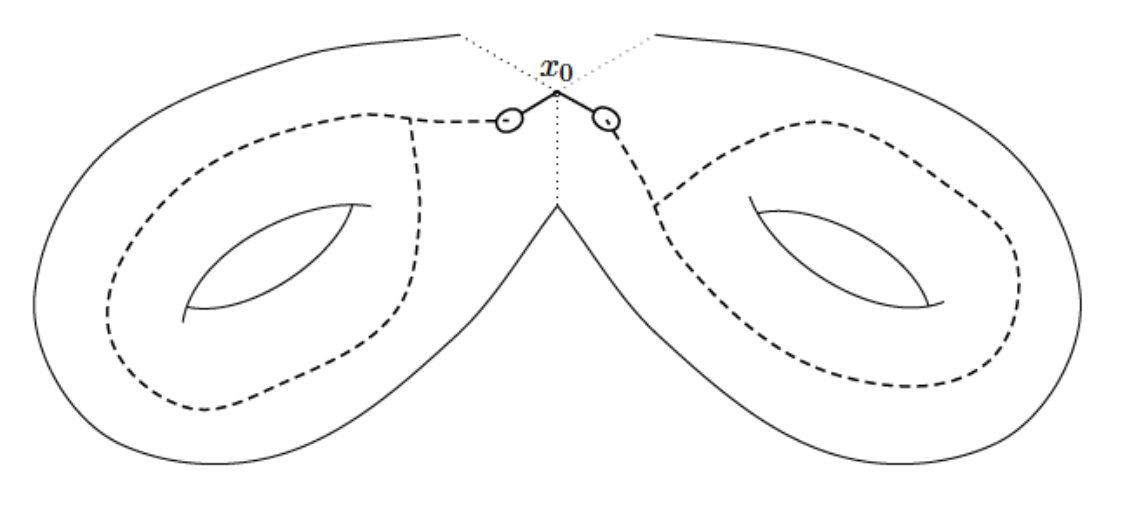}}%
  \subfigure[Linking in Each Handle]%
    {%
    \label{SizeR2}
    \includegraphics[width=0.45\textwidth]{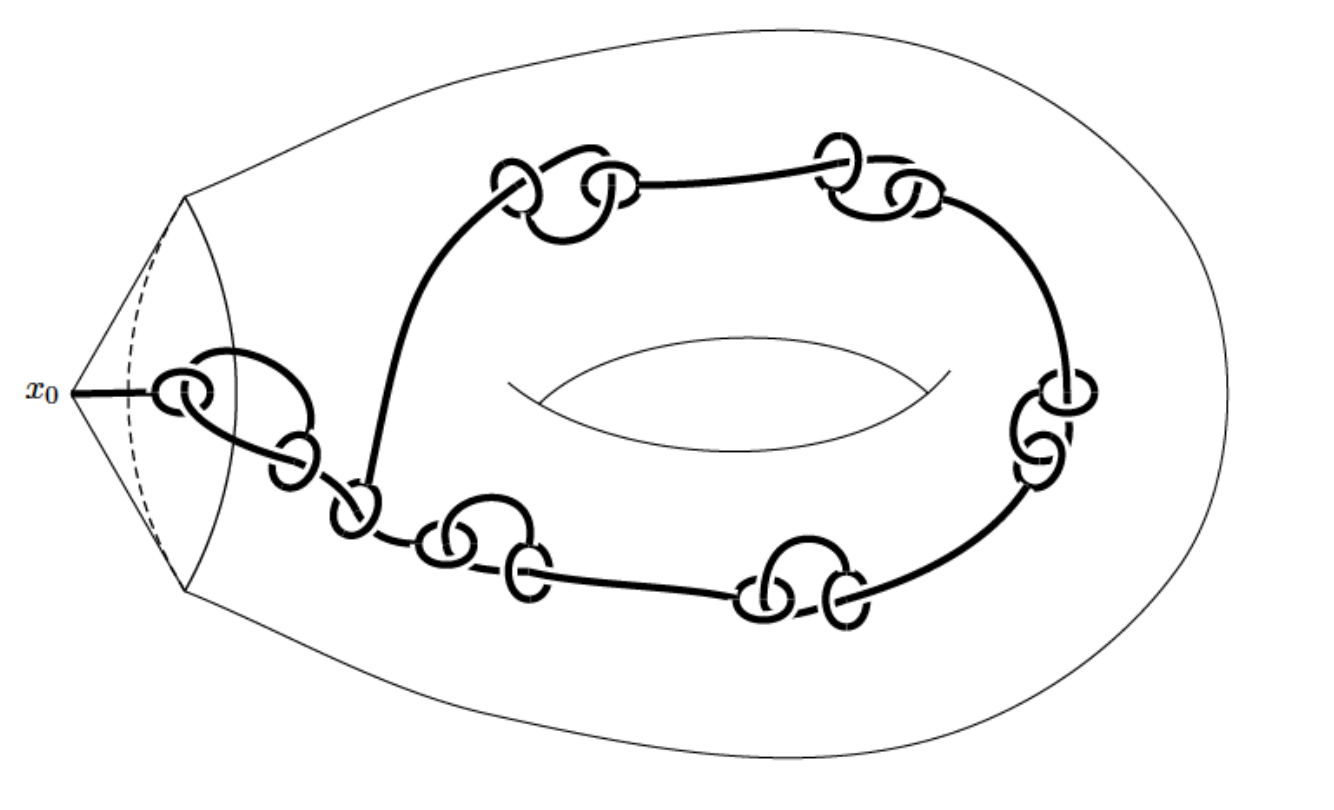}
    }\\%
  \end{center}
  \caption{%
    Size Replacement}%
  \label{SizeR}
\end{figure}

(See Figure \ref{SizeR}).
Let $H$ be an unknotted genus $k$ solid handlebody in $S^3$. View $H$ as a solid ball with $k$ 1-handles attached as in Figure \ref{Replacement1}. Replace $H$ by a single smaller genus $k$ handlebody $H^{\prime}$ contained in the solid ball (see Figure \ref{SizeR1}), and by  chains of smaller unknotted genus 2 handlebodies in each handle of $H$ as in Figure \ref{SizeR2}. This replacement can be done so the diameter of each of the new handlebodies in $H$ is less than half the diameter of $H$. 

\begin{remark} 
\label{3cellRemark}
Let $B$ be a 3-cell containing $H$, and let $H_1, H_2, \ldots H_k$ be the new handlebodies in $H$. Then there are pairwise disjoint 3-cells $B_1, B_2, \ldots B_k$ in $B$ so that $H_i\subset B_i$.  This follows from the fact that the genus 2 handlebodies around each handle are not linked in $B$ even though they are linked in $H$. See \cite{Sk86} and \cite{GRZ06} for more details on this.
\end{remark}
 \textbf{Note:} The inductive hypotheses imply that there is a nested sequence of unions of pairwise disjoint 3-cells associated with the Cantor set constructed in the next section. This collection of 3-cells however does not imply that the constructed Cantor set is tame, since the diameters of the 3-cells do not go to 0. 

\section{The Construction}
\label{constructionsec}

Fix a sequence $\mathfrak{S}=(n_1,n_2, \ldots)$ where each $n_i$ is either an integer greater than 1 or is $\infty$. 
We will construct inductively a defining sequence $S_1$, $S_2$, $\ldots\ $ for a Cantor set $X=X(\mathfrak{S})$ in $S^3$.  Equivalently, we specify an increasing union of open $3-$manifolds $N_i$ where each $N_i=S^3\setminus S_i$. At the same time, we associate with each component $M_{(i,j)}$ of $S_i$  the specific term $n_j$ of the sequence $\mathfrak{S}$ above.

\textbf{Inductive hypotheses:}  
\begin{enumerate}
\item[\textbf{IH 1:}] The components $M_{(i,1)},\ldots M_{(i,m(i))}$ of $S_i$ are unknotted handlebodies of genus 2 or greater which are contained in pairwise disjoint 3-cells in $S^3$.
 
\item[\textbf{IH 2:}] The genus of $M_{(i,j)}$ is less than or equal to $n_j$.
\end{enumerate}

The construction will show that the components of $S_{2k+1}$ are obtained from $S_{2k}$ by
a suitable replacement construction performed on the components of $S_{2k}$. 
The components of
$S_{2k}$ will be obtained by replacing each component of 
$S_{2k-1}$ by an appropriate chain of linked handlebodies. All except
possibly one handlebody in this chain will have genus 2.

\textbf{Stage 1:} To begin the construction, let $S_1$ consist of a single component $M_{(1,1)}$, where $M_{(1,1)}$ is an unknotted genus $2$ handlebody in $S^3$. Note that $n_1$ is greater than or equal to the genus of $M_{(1,1)}$. The inductive hypotheses are clearly satisfied  since $n_1\geq 2$.

\textbf{Stage $\mathbf{k+1}$ if $\mathbf{k}$ is odd:}

By  the inductive hypothesis, every component of
    $S_k$ is an unknotted handlebody of genus  2 or greater. Let $M_{(k,i)}$ be a component of $S_k$. Again, by the inductive hypothesis, the genus of $M_{(k,i)}$ is less than or equal to $n_i$. If the genus of $M_{(k,i)}$ is less than $n_{i}$, perform a genus replacement on $M_{(k,i)}$. Note that this replaces a component $M_{(k,i)}$ at stage $k$ by a component $M_{(k+1,i)}$ at stage $k+1$. This genus of $M_{(k+1,i)}$ is then  (genus of $M_{(k,i)}$)$+1$.
So the genus of $M_{(k+1,i)}$ is less than or equal to $n_i$.

If the genus of $M_{(k,i)}$ is equal to $n_i$, let $M_{(k+1,i)}$ be $M_{(k,i)}$ with a small open neighborhood of the boundary removed.  Since the genus of $M_{(k+1,i)}$ is the genus of $M_{(k,i)}$, this genus is still less than or equal to $n_i$. In either case, inductive hypothesis \textbf{IH 2} is satisfied at stage $k+1$.

By construction, each component $M_{(k+1,i)}$ is an unknotted handlebody of genus at least 2. Since each new component is in the interior of a component from $S_n$, the inductive hypothesis \textbf{IH 1} at stage $k$ guarantees that the inductive hypothesis \textbf{IH 1} is still satisfied at stage $k+1$.

\textbf{Stage $\mathbf{k+1}$ if $\mathbf{k}$ is even:}

By  the inductive hypothesis, every component of
$S_k$ is an unknotted  handlebody of genus  2 or greater. Let $M_{(k,i)}$ be a component of $S_k$ of genus $g$. Again, by the inductive hypothesis,  $g$ is less than or equal to $n_i$. Perform a size replacement on $M_{(k,i)}$, replacing $M_{(k,i)}$ by a single genus $k$ 
smaller copy  (labelled $M_{(k+1,i)}$) of $M_{(k,i)}$, and by 6g unknotted genus 2 handlebodies as in  Figure \ref{SizeR}. 

This procedure performed on each component of $S_k$ produces smaller copies of $M_{(k,i)}, 1\leq i\leq m(n)$, labeled $M_{(n+1,i)}$, and a collection of unlabeled genus 2 handlebodies. Arbitrarily label these new genus 2 handlebodies as
$M_{(k+1,m(k)+1)}\ldots M_{(n+1,m(k+1))}$.

By \textbf{IH 1} at stage $k$, by the construction above, and by Remark \ref{3cellRemark}, 
\textbf{IH 1} is satisfied at stage $k+1$. 

If $1\leq i\leq m(k)$, the genus of $M_{(k+1,i)}$ is equal to the genus of $M_{(k,i)}$ which is less than or equal to $n_i$, by \textbf{IH 2} at stage $k$. If $i\geq m(k)+1$,  the genus of $M_{(k+1,i)}$ is equal to 2 which is less than or equal to $n_i$ by the assumption about the terms of $\mathfrak{S}$.
In either case, \textbf{IH 2} is satisfied at stage $k+1$.

\begin{figure}[hbt]
\begin{center}
\includegraphics[width=.5\textwidth]{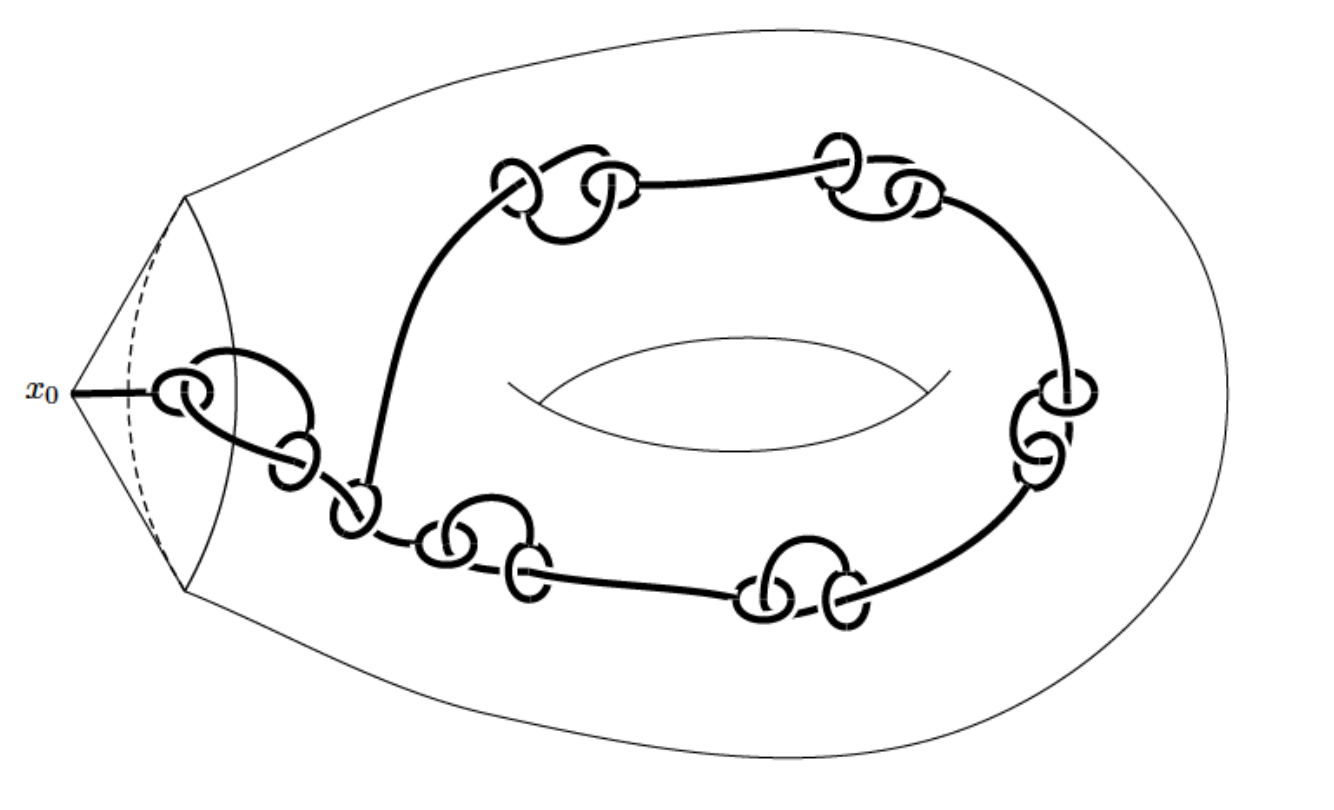}
\end{center}
\caption{Linking along the spine of some handle of $N$}
\label{Linking}
\end{figure}

This completes the inductive description of the defining
sequence. Define the Cantor set $C=C_{\mathfrak{S}}$ associated with the sequence $\mathfrak{S}$ to be $C=C_{\mathfrak{S}}=\bigcap_i M_i\}.$

In the next section, we will show that   $C=C_{\mathfrak{S}}$ is the required Cantor set needed for proving Theorem \ref{CantorTheorem}.

\section{Proof of the Main Theorems}
\label{proofsec}
For a sequence $\mathfrak{S}=(n_1,n_2, \ldots)$ where each $n_i$ is either an integer greater than 1 or is $\infty$, let $C=C_{\mathfrak{S}}$ be the Cantor set constructed in the previous section.
Note that $C$ is indeed a Cantor set since each component of $M_n$ contains more than one component of $M_{n+2}$ and since the diameter of the components goes to 0 as $n$ goes to $\infty$.
We now verify that $C$ satisfies the conditions listed in Theorem \ref{CantorTheorem}.

\subsection{
Simple connectivity of the complement}
Let $\gamma\colon S^1\to S^3\setminus C$. The set $\gamma(S^1)$
is compact and misses $C$ so there exists an even $n$ large enough such
that $\gamma(S^1)\cap M_n=\emptyset$. 

By \textbf{IH 1}, the components $M_{(n,1)},\ldots M_{(n,m(n))}$ of $M_n$ are unknotted handlebodies of genus 2 or greater which are contained in pairwise disjoint 3-cells in $S^3$. Since the components are cubes with unknotted handles and lie in disjoint 3-cells,
the fundamental group of the complement of the components is
generated by the meridional curves on the components. It
therefore suffices to show how any meridional loop (say $J$) of
some component $N$ can be shrunk to a point in the complement of
the components at the next stage contained in $N$.

\begin{figure}[htb]
\begin{center}
\includegraphics[width=.5\textwidth]{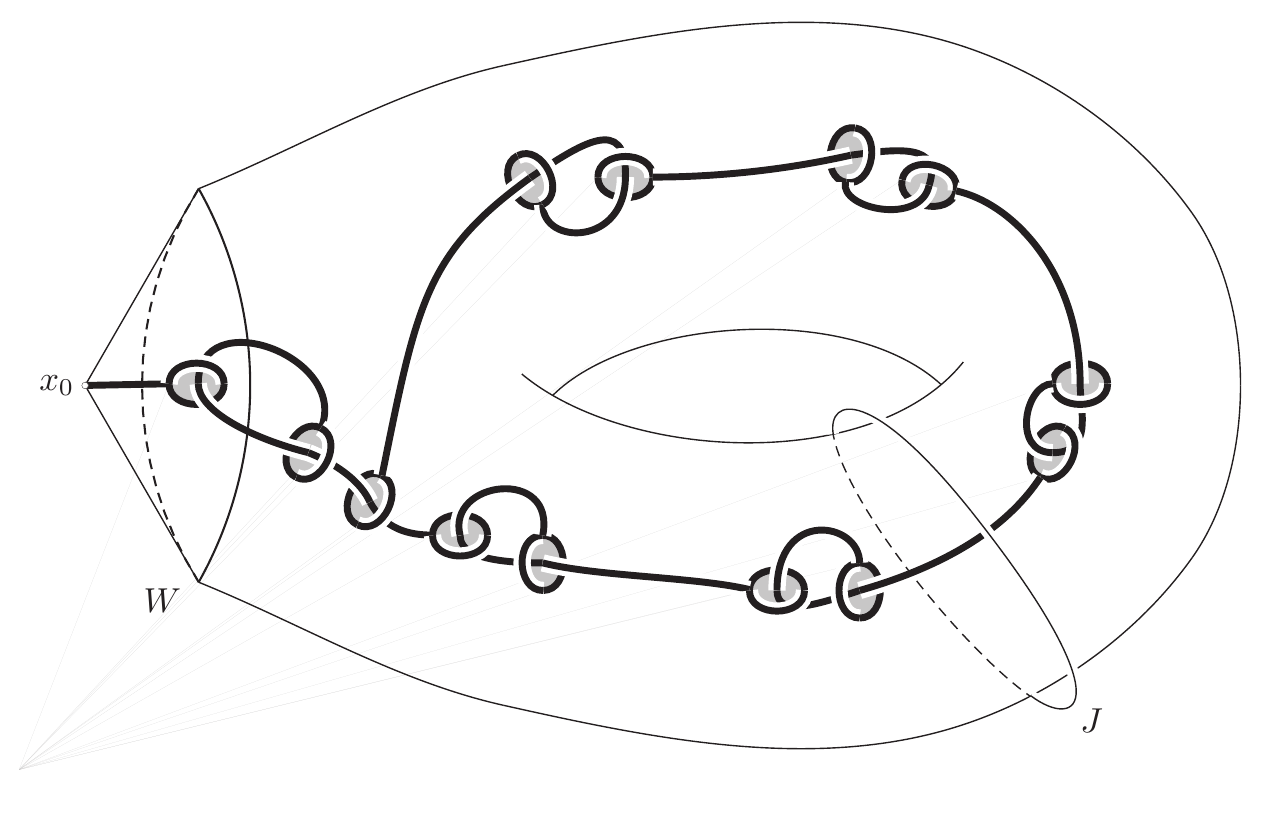}
\end{center}
\caption{Simple Connectivity of the Complement}
\label{AddDiscs}
\end{figure}

By construction it is clear that $J$ can be moved
in $N\setminus M_{n+1}$ to the loop $W$ on $\text{Bd}N$ (see Figure
\ref{AddDiscs}) and then moved off $N$. Hence
$[J]=0\in\pi_1(S^3\setminus C)$.\qed

\subsection{The countable dense subset}   Fix a point $x\in C$ and let $M(j,x(j))$ be the component of $M_j$ containing $x$. Then
$$
 M(1,x(1))\supset M(2,x(2))\supset \ldots \supset M(n,x(n))\supset \ldots,
\text{ and } 
$$
$$
x=\bigcap_n M(n,x(n)).
$$

  {For a fixed $i>0$, define the point $x_i\in C$ as follows. Choose $n$ so that  the number $m(n)$ of components in $M_n$, is greater than $i$. Then $M(k,i)$ is a component 
of $M_k$ for each $k\geq n$ and by construction  $M(k+1,i)\subset M(k,i)$.  Let $x_i=\displaystyle\bigcap_{k\geq n} 
M(k,i)$. The countable dense subset of $C$ that we are looking for is the subset $D=\{x_1,x_2.\ldots\}$. The set $D$ is dense in $C$ since each component $M(k,j)$ of $M_k$ contains a point of $D$, namely $x_j$.

\subsection{Genus at $x_j$ when $n_j$ is finite.} 
 {Note that $x_j=\displaystyle\bigcap_{k\geq n} 
M(k,j)$ where $n$ is chosen so that $m(n)\geq j$. 
The genus of $M(n,j)$ is less than or equal to $n_j$, and the genus modification at subsequent stages guarantees that there is an $n^{\prime}\geq n$ so that the genus of $M(n^{\prime},j)$ is equal to $ n_j$. Then for each $\ell>n^{\prime}$, the genus of $M(\ell,j)$ is also equal to $n_j$. It follows from the definition of local genus that
the local genus $g_{x_j}(C)$ of $C$ at $x_j$  is then less than or equal to $n_j$.

 It remains to be shown that $g_{x_j}(C)$ is greater than or equal to $n_j$. 

Choose the component  $N=M(2i+1, j)$ of $M_{2i+1}$ containing $x_j$, where $m(2i+1)>j$ and $2i+1$ is large enough so that the genus of $N$ is equal to $n_j$.  Then $N$ is a union of genus 1 handlebodies as in the previous section. Let T be one of these genus 1 handlebodies. By construction we have that $\text{Bd}(T )\cap C$ is the singleton $x_j$. Let W be a loop in Bd(T) that bounds a disc in Bd(T) containing $x_j$ in its interior as in Figure \ref{LoopW}.

\begin{figure}[ht]
\begin{center}
  \subfigure[The loop $W$]%
    {%
    \label{LoopW}
    \includegraphics[width=0.45\textwidth]{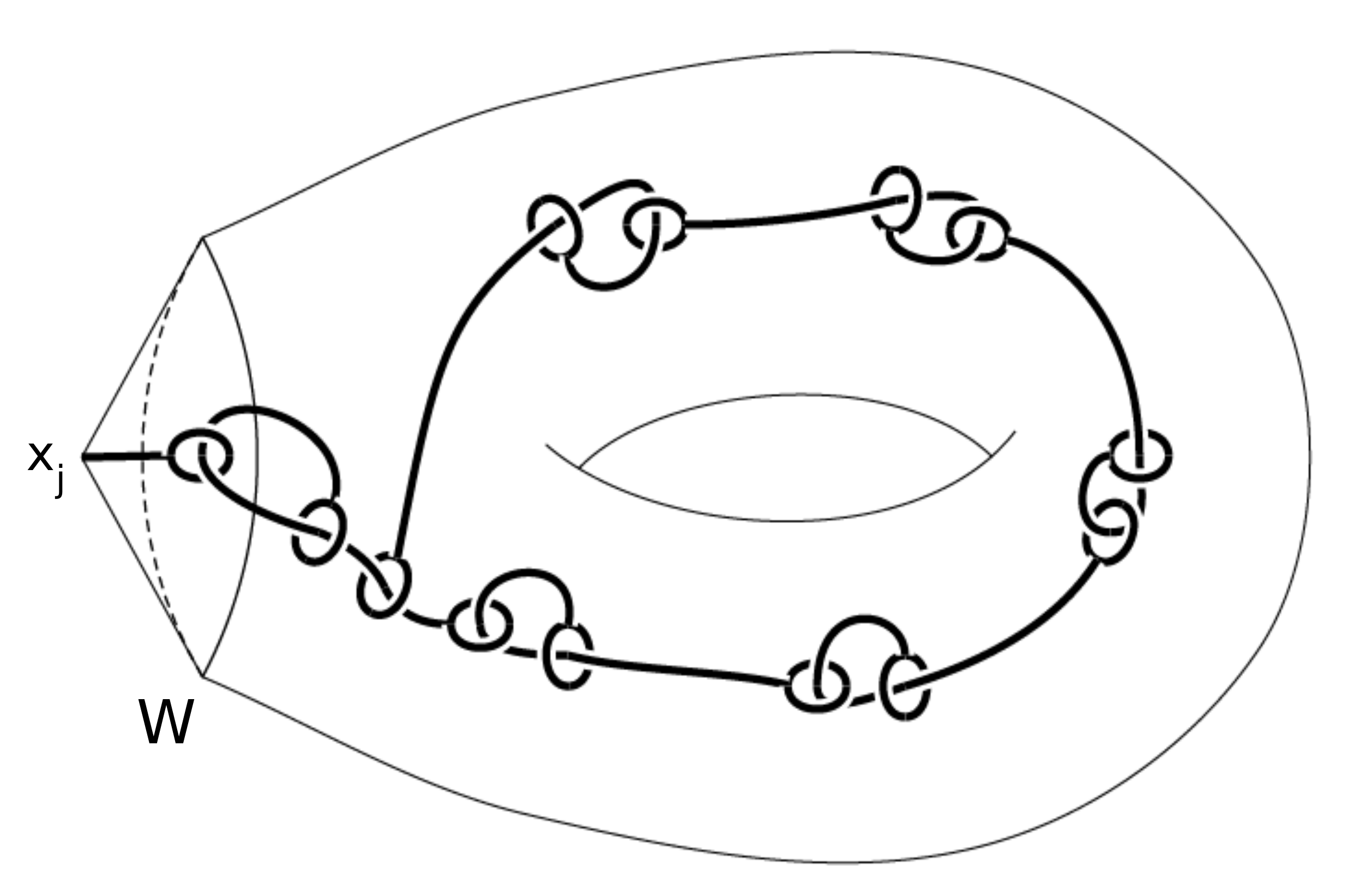}
    }%
  \subfigure[$T$ with filled in disks]%
    {%
    \label{Fill}
    \includegraphics[width=0.5\textwidth]{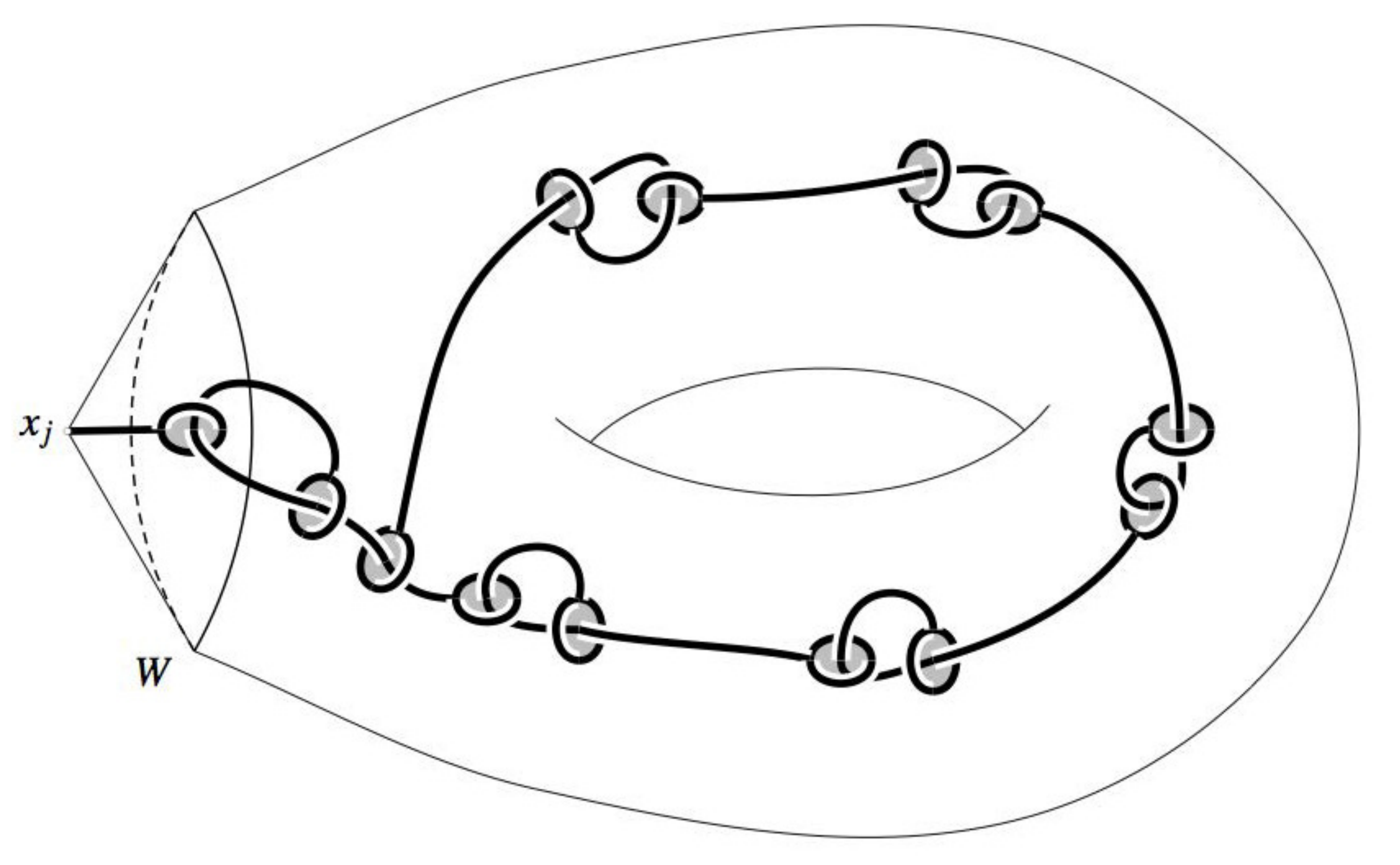}
    }\\%
  \end{center}
  \caption{%
    The handlebody $T$}%
  \label{HandlebodyT}
\end{figure}

By Theorem \ref{Slicing}, to show that the local genus of $x_j$ in $C$ is at least $n_j$, it suffices to show that the local genus of $x_j$ in $C\cap T$ is at least 1. For this, we will need two technical lemmas.

\begin{lem}
\label{SpineLemma}
If there exists a 2-disc $D\subset T$ such that $D\cap
M_{r+1}=\emptyset$ for some $r>2i+1$, and $Bd(D)=W$, then there exists a 2-disc
$D\thinspace'\subset T$ such that $Bd(D\thinspace')=Bd(D)$ and
$D\thinspace'\cap M_{r}=\emptyset$.

\end{lem}
\textbf{Proof:} This is essentially Lemma 5.1 in \cite{GRZ06}. Complete details of the proof are provided there.  The main idea is that if $D\cap
M_{r+1}=\emptyset$, $D$ can be adjusted  so as to miss a core of each component of $M_r\cap T$. This allows formation of a new disc $D^{\prime}$ with the same boundary so that $D\thinspace'\cap M_{r}=\emptyset$. \qed

\begin{lem}
\label{GenusZeroLemma}
Let $W$ be a loop on $Bd(T)$
as in Figure \ref{LoopW} and let $x_j$ be the point in $C\cap
Bd(T)$.  If $g_{x_j}(C\cap T)=0$, then $W$ bounds a disc $D$ in
$T$ missing $C$.
\end{lem}
\textbf{Proof:} This is essentially Lemma 5.2 in \cite{GRZ06}. Complete details of the proof are provided there.  The argument there is for a component $N$ of genus greater than or equal to 3. Since the argument is actually applied separately to each 1-handle of $N$, the same argument works here where the genus of $N$ is greater than or equal to 2. \qed

Now assume that $g_{x_j}(C\cap T)=0$. Then by Lemma \ref{GenusZeroLemma},
the loop $W$ bounds a disc $D$ missing $C$, and so missing $M_r\cap T$ for some large $r$. Repeated application of Lemma \ref{SpineLemma} then implies that $W$ bounds a disc $D\thinspace^{\prime}$ so that $D\thinspace'\cap M_{2i+2}=\emptyset$. An argument similar to that used in Section 6.2 of \cite{GRZ06} then shows that $W$ bounds a disc also missing the filled in disks in Figure \ref{Fill}. This disk can then be pushed to the boundary of $T$ minus $x_j$. Since $W$ is nontrivial in the boundary of $T$ minus $x_j$, this cannot happen. It follows that $g_{x_j}(C\cap T)\geq 1$}.

This completes the argument that the local genus of $C$ at $x_j$ is equal to $n_j$ when $n_j$ is finite.}\qed

\subsection{Genus at $x_j$ when $n_j$ is infinite.}  It suffices to show that $g_{x_j}(C)$ is greater than or equal to $K$ for each positive integer $K$.

Fix a positive integer $K$ and choose a stage $M_n$ of the defining sequence for $C$ so that  $m(n)>j$ and so that the genus $g$ of $N=M(n,j)$, is at least $K$. Then $N$ is a union of $g$ genus 1 handlebodies. Let T be one of these genus 1 handlebodies. By construction, even though the genus of the handlebody containing $x_j$ in $M_r$ increases as $r$ goes to infinity,  we still have that $\text{Bd}(T )\cap C$ is the singleton $x_j$.

Exactly the same argument as above (in the case when $n_j$ was finite) can now be applied to conclude that the local genus of $x_j$ in $C\cap T$ is at least one. It follows that the local genus of $x_j$ in $C\cap N$ is at least $K$ which is what we needed to show. So the local genus of $x_j$ in $C$ is infinity.\qed

\subsection{Genus at points not in the dense set.
}
If $x\in C\setminus D$, then for each $N\in \mathbb{Z}_+$, there is an $N^{\prime}\geq N$ so that the genus of $M(N^{\prime},x(N^{\prime}))$ is two. To see this, consider the component $M(N,i)$ of $M_N$ that contains $x$. If the genus of $M(N,i)$ is two, we are done. If not, there must be a first stage $N+k>N$ such that $x(N+k)\neq i$. If there is no such stage, then $x=x_i\in D$ which can't happen. For this stage $N+k$, $M(N+k, x(N+k))$ has genus two by the construction of $C$. It follows that in the sequence whose intersection is $x$,
$$
 M(1,x(1))\supset M(2,x(2))\supset \ldots \supset M(n,x(n))\supset \ldots\ ,
$$
infinitely many of the terms have genus two. It then follows from the definition of local genus that
$g_x(C)\leq 2$.\qed

\begin{remark}
\label{genusnotzero}
An inductive argument using techniques similar to those of Lemma \ref{SpineLemma}
can be used to show the following. Let  $M(k,i)$ be a component of $M_k$. Then there are no 2-spheres missing $C$ in the interior of 
$M(k,i)$  that separate $M(k,i)\cap C$. This in turn can be used to show that the local genus of each point in $C$ is strictly greater than 0. So the local genus at each point of $C\setminus D$ is either 1 or 2.
\end{remark}

This completes the proof of Theorem \ref{CantorTheorem}. By Lemma \ref{equivalence}, Theorem \ref{EndTheorem} (the main theorem) follows as an immediate corollary.\qed

\section{Questions}
\label{questionsec}

Both Theorems \ref{EndTheorem} and \ref{CantorTheorem} contain the requirement that each term in the sequence $\mathcal{S}$ is greater than one. Additionally, the genus of each end not in the dense set $D$ (Theorem \ref{EndTheorem})  and of each point of $C$ not in $D$ (Theorem \ref{CantorTheorem}) are shown to be less than two. This leads to two interesting questions, listed below.

\begin{que}

Can the restrictions  that $n_i\geq 2$ in Theorem \ref{EndTheorem} and Theorem \ref{CantorTheorem}  be removed? That is, is it possible to construct simply connected examples with certain points in the dense set specified to have genus 1?
\end{que}
\begin{que}
By Remark \ref{genusnotzero}, the genus of points not in  the dense sets from  
Theorem \ref{EndTheorem} and Theorem \ref{CantorTheorem} is either 1 or 2. Can it be shown that their genus  is exactly equal to 2?
\end{que}
It seems that new construction techniques need to be developed to answer the first question since the specific genus 2 handlebodies used in the size replacement modification are essential to the argument that the complement of the Cantor set constructed is simply connected.

\end{document}